\documentclass[12pt,reqno]{amsart}

\usepackage%[dvips] %including this option causes pdflatex to fail to
{graphicx} 

\usepackage{amssymb,  %this is needed for \measuredhornangle
hyperref
}

\usepackage{breakurl}   %this allows for url breaks in pdf

\usepackage[export]{adjustbox}  %this is necessary for graphics, figures

\usepackage{epigraph}

\usepackage{setspace}  %this is needed for double-spacing

\usepackage{apalike} %to change citations to parentheses instead of
\theoremstyle{definition}

\numberwithin{equation}{section}

\title{Episodes from the history of infinitesimals}

\author[Mikhail Katz]{Mikhail G. Katz} \address{M. Katz, Department of
  Mathematics, Bar Ilan University, Ramat Gan 5290002 Israel}
\email{katzmik@math.biu.ac.il}


\begin{document}

\begin{abstract}
Infinitesimals have seen ups and downs in their tumultuous history.
In the 18th century, d'Alembert set the tone by describing
infinitesimals as chimeras.  Some adversaries of infinitesimals,
including Moigno and Connes, picked up on the term.  We highlight the
work of Cauchy, No\"el, Poisson and Riemann.  We also chronicle
reactions by Moigno, Lamarle and Cantor, and signal the start of a
revival with Peano.
\end{abstract}

%\doublespacing

\thispagestyle{empty}

%\huge


\maketitle
%\tableofcontents

%\today

\section{Chimeras from Moigno to Connes}  
\label{s2}

The history of infinitesimals involved a number of philosophical and
religious controversies already during the 17th century and even
earlier.
%
%this should include a footnote on jesuat.tex
%
Our intention here is to deal with the post-Leibnizian period, not
including Leibniz.  The reason is that Leibniz has been amply covered
elsewhere, including the most recent exchange between Arthur and
Rabouin \cite{Ar24} and Katz and Kuhlemann \cite{23f}.  For older
studies, see Goldenbaum and Jesseph \cite{Go08} and work cited
therein.

In the 18th century, d'Alembert used the deprecating epithet
\emph{chimera} (referring to an illusory thing) for infinitesimals.%
\footnote{See Lamand\'e \cite[52]{La19}.}
The epithet found an eager following in the ensuing centuries.  In
19th century France, infinitesimals were panned as \emph{chimeras} by
jesuit-trained Abb\'e Moigno.  Moigno wrote:%
\footnote{\label{f4}A measure of Moigno's mathematical accomplishment
can be gleaned from an obituary in \emph{Nature}: ``In addition to
theology he studied with great enthusiasm both the physical and
mathematical sciences; in these he made rapid progress and in 1828
arrived at a new mode of getting the equation to [sic] the tangent
plane to a surface'' \emph{Nature} \cite[291]{Mo84}.}
\begin{enumerate}\item[]
In effect, either these magnitudes, smaller than any given magnitude,
still have substance and are divisible, or they are simple and
indivisible: in the first case their existence is a \emph{chimera},
since, necessarily greater than their half, their quarter, etc., they
are not actually less than any given magnitude; in the second
hypothesis, they are no longer mathematical magnitudes, but take on
this quality, this would renounce the idea of the continuum divisible
to infinity, a necessary and fundamental point of departure of all the
mathematical sciences.%
\footnote{Moigno as translated by Schubring \cite[456]{Sc05}; emphasis
added.}
\end{enumerate}
In the 20th century, the term was picked up by Moigno's compatriot
Alain Connes:
\begin{enumerate}\item[]
A nonstandard number is some sort of \emph{chimera} which is
impossible to grasp and certainly not a concrete object.  In fact when
you look at nonstandard analysis you find out that except for the use
of ultraproducts, which is very efficient, it just shifts the order in
logic by one step; it's not doing much more.%
\footnote{Connes \cite[14]{Co00}; emphasis added.}
\end{enumerate}
Connes repeatedly voiced opinions about infinitesimals in the course
of three decades.  His views are analyzed in Kanovei et
al.~\cite{13c}, Katz and Leichtnam \cite{13d}, and Sanders \cite{18l}.

\section{From Cauchy to Riemann}  
\label{s2b}

Whether `hallucinations' (Mario Bettini)%
\footnote{See Alexander \cite[159]{Al15}.}
or `chimeras' (d'Alem\-bert), infinitesimals in the hands of Leibniz,%
\footnote{See Katz and Sherry \cite{13f}, Katz et al.~\cite{23h}, @@.}
Bernoulli, Euler,%
\footnote{For an analysis of Euler's use of infinitesimals, see Bair
et al.~\cite{17b}.}
and others fared fairly well until and including the first half of the
19th century, when leading mathematicians like Cauchy, Poisson, and
Riemann still made routine use of infinitesimal techniques.%
\footnote{For a summary of the debate over this historical period, see
%\cite{23a}, 
Bair et al.~\cite[Section~2.3]{22a}.}

\subsection{Cauchy's `epsilontic' epsilon}
\label{s31}

We explore the nature of epsilons%
\footnote{The plural is significant; see note~\ref{f6}.  See further
in Bair et al.~\cite[Section~2]{22a}.}
occurring in Cauchy's work.  A prototype of an~$\varepsilon,\delta$
argument occurs in his 1821 textbook \emph{Cours d'Analyse}.%
\footnote{Cauchy \cite{Ca21}, Section~2.3, Theorem~1.}
Here Cauchy assumes that the difference~$f(x+1)-f(x)$ converges
towards a limit~$k$ for increasing values of~$x$.  He argues that the
ratio~$\frac{f(x)}{x}$ converges toward the same limit.%
\footnote{See Bradley and Sandifer (transl.) \cite[35]{BS}.}
He chooses~$\varepsilon>0$ and notes that we can give the number~$h$ a
value large enough so that, when~$x\geq h$, the
difference~$f(x+1)-f(x)$ is always contained between~$k-\varepsilon$
and~\mbox{$k+\varepsilon$}.  He then arrives at the formula
\[
\frac{f(h+n)-f(h)}{n}=k+\alpha,
\]
where~$\alpha$ is a quantity contained between~$-\varepsilon$
and~$+\varepsilon$.  He concludes that the ratio~$\frac{f(x)}{x}$ has
for its limit a quantity contained between~$k-\varepsilon$
and~$k+\varepsilon$.  Here~$\varepsilon$ denotes an ``ordinary''
number (not infinitesimal),%
\footnote{Grabiner \cite[8]{Gr81}; Bair et
al.~\cite[Section\;3.8]{19a}.}
as it does in a 1823 argument in the same vein, starting with the
words
\begin{quote}
``Denote by~$\delta$,~$\varepsilon$ two very small numbers.''%
\footnote{Cauchy as translated by Cates in \cite[33]{Ca19}.}
\end{quote}
Neither of these arguments involves infinitesimals.

\subsection{Cauchy's infinitesimal epsilon}
\label{s31b}

In 1826, Cauchy studies the notions of the center and radius of
curvature of curves.%
\footnote{See Bair et al.~\cite{20b}, Section~5, pp.\,133--136.}
He defines the \emph{angle de contingence}~$\pm\Delta\tau$ as the
angle between the two tangent lines of an arc~$\pm\Delta s$ at its
extremities.  The term \emph{angle of contingence} (sometimes referred
to as a \emph{hornangle} or an \emph{angle of contact}) has
traditionally referred to a curvilinear angle incomparably smaller
than ``ordinary'' rectilinear angles, and therefore indicative of
non-Archimedean behavior.  The matter was debated extensively in the
16--17th centuries and even earlier; Leibniz used the term in a
similar non-Archimedean sense.%
\footnote{For details on Leibniz's work on hornangles see Katz et
al.~\cite[Section~4.5]{24c}.}

Cauchy considers the normals to the curve at the extremities of an arc
starting at the point~$(x,y)$.  He gives \emph{two} definitions of
both the center of curvature and the radius of curvature:
\begin{enumerate}\item[]
The distance from the point~$(x,y)$ to the intersection point of two
normals is appreciably equivalent to the radius of a circle which
would have the same curvature as the curve.%
\footnote{``La distance du point~$(x, y)$ au point de rencontre des
deux normales est sensiblement \'equivalente au rayon d'un cercle qui
aurait la m\^eme courbure que la courbe'' Cauchy \cite[98]{Ca26}.}
\end{enumerate}
The intersection of the two normal lines produces a point which will
generate the \emph{center} of curvature; the distance between~$(x,y)$
and the center is the \emph{radius} of curvature.  In order for the
intersection of the two normals to generate the center of curvature,
the distance~$\Delta s$ between them must be infinitesimal, similarly
to~$\Delta\tau$.

Next, Cauchy chooses an infinitesimal number~$\varepsilon$ and
exploits the law of sines to write down a relation that will give an
expression for the radius of curvature:
\begin{enumerate}\item[]
si l'on d\'esigne par~$\varepsilon$ un \emph{nombre} infiniment petit,
on aura%
\footnote{Ibid.; emphasis added.}
\begin{equation}
\label{e61b}
\frac{\sin\left(\frac{\pi}{2}\pm\varepsilon\right)}{r}=
\frac{\sin(\pm\Delta\tau)}{\sqrt{\Delta x^2+\Delta y^2}}
\end{equation}
\end{enumerate}
Note that here Cauchy explicitly describes his~$\varepsilon$ as an
infinitely small number.  He then passes to the limit to obtain
\begin{equation}
\label{e62b}
\frac{1}{\rho} = \pm \frac{d\tau}{\sqrt{dx^2+dy^2}}
\end{equation}
(ibid., p.\;99).  To pass from formula~\eqref{e61b} to
formula~\eqref{e62b}, Cauchy replaces the infinitesimals~$\Delta
x$,~$\Delta y$, and~$\sin\Delta\tau$ by the corresponding
differentials~$dx$,~$dy$, and~$d\tau$.%
\footnote{Note that Cauchy's differentials were arbitrary numbers (not
necessarily infinitesimals).  In modern infinitesimal analysis, to
calculate the derivative of~$y=f(x)$, one would choose a nonzero
infinitesimal increment~$\Delta x$, compute the corresponding
change~$\Delta y = f(x+\Delta x) - f(x)$, and form the
quotient~$\frac{\Delta y}{\Delta x}$.  The
derivative~$f'(x)=\frac{dy}{dx}$ is then computed as the standard part
of the quotient~$\frac{\Delta y}{\Delta x}$.  Here the infinitesimal
differentials~$dx$ and~$dy$ in the quotient~$\frac{dy}{dx}$ are chosen
as follows:~$dx=\Delta x$ whereas~$dy=f'(x)dx$.  Of course, the
symbols~$\Delta x$ and~$\Delta y$ are also used in traditional
non-infinitesimal calculus, where they have a different meaning.}
Here the expression~$\sin(\frac{\pi}{2}\pm\varepsilon)$ is infinitely
close to~$1$ whereas~$r$ is infinitely close to the radius of
curvature ~$\rho$, justifying the replacements in the left-hand side
of Cauchy's equation.  Cauchy's~$\Delta x$,~$\Delta y$,~$\Delta
s$,~$\Delta\tau$, and~$\varepsilon$ are all genuine infinitesimals.%
\footnote{\label{f6}Cauchy's~$\varepsilon$ here is unlike the
symbol~$\varepsilon$ exploited in \cite{Ca21}, Section~2.3, Theorem~1
and in 1823, as discussed in Section~\ref{s31}.}

Cauchy's 1853 sum theorem for convergent series of continuous
functions is another illustration of his use of genuine
infinitesimals.%
\footnote{See Bascelli et al.~\cite{18e}.}
The tendency of some historians to deny that Cauchy used genuine
infinitesimals is analyzed in Bair et al.~\cite{17a} and earlier in
Spalt \cite{Sp81}.  The rival interpretations are summarized
in~Katz \cite{21f}.%
\footnote{See Bair et al.~\cite[Section~3.1]{22a} for a summary of the
debate over the sum theorem.}
Other contemporary advocates of infinitesimals are discussed in
Section~\ref{s33b}.

\subsection{Poisson, No\"el, Manilius}
\label{s33b}

A pedagogical dispute in Belgium and Luxembourg in the middle of the
19th century concerned the teaching of the calculus.  No\"el and
others advocated the use of infinitesimals.  No\"el's position on
infinitesimals was consistent with that of Poisson.  Indeed,
\begin{enumerate}\item[]
Poisson {\ldots}\;was able to promote them to far-reaching
dissemination and impact because of his central position within the
French educational system.%
\footnote{Schubring \cite[575]{Sc05}.}
\end{enumerate}
Another leading mathematician who endorsed genuine infinitesimals
during this period is Coriolis.%
\footnote{See Grattan-Guinness \cite[1262]{Gr90}.}

Lamarle and others opposed the approach using infinitesimals (see
Section~\ref{s33}).  The dispute is part of a 19th century
infinitesimal lore%
\footnote{For a discussion of 19th century infinitesimal lores, see
Bair et al.~\cite[Section~2.7]{22a}.}
that featured markedly different conceptions of infinitesimals,
including No\"el's conception thereof as inassignable magnitudes
violating what later became known as the Archimedean property.

We will start with No\"el and Manilius, and analyze the views of
Lamarle and Paque in Section~\ref{s33}.  Jean-Nicolas No\"el
(1783--1867) defined infinitesimals as follows in 1855:
\begin{enumerate}\item[] 
A number is said to be infinitely small if it is smaller than the
least assignable part of the unity.  Such a number is therefore
absolutely inappreciable due to its smallness and could never be
expressed in digits.%
\footnote{``Un nombre est dit infiniment petit lorsqu'il est moindre
que la plus petite partie assignable de l'unit\'e. Un tel nombre est
donc absolument inappr\'eciable par sa petitesse et ne pourra jamais
s'exprimer en chiffres'' No\"el \cite[28]{No55}.  See also Bockstaele
\cite[13]{Bo66} and Bair and Mawhin~\cite[48]{Ba19}.}
\end{enumerate}
In an earlier text, No\"el had explicitly mentioned the
assignable \emph{vs} inassignable dichotomy.%
\footnote{``D’apr\`es la v\'eritable d\'efinition du rapport, deux
grandeurs continues de m\^eme nature ont toujours un rapport unique,
exprimable ou inexprimable (rationnel ou irrationnel); elles ont donc
aussi toujours une mesure commune, \emph{assignable ou
inassign\-able}, finie ou infiniment petite, et toujours inconnue,
dans ce dernier cas'' (No\"el as quoted by Mawhin in \cite{Ma17};
emphasis added).}
No\"el's assignable/inassignable distinction is mentioned by
Bockstaele.%
\footnote{See Bockstaele \cite[11]{Bo66}.}
Both No\"el's dichotomy of assignable \emph{vs} inassignable
magnitudes and his definition of infinitesimal as smaller than every
assignable magnitude are essentially the Leibnizian concepts; see
Section~\ref{f9}.%
\footnote{See further in Bair et al.~\cite[Section 2.7]{22a}.}

No\"el's ally Jean Joseph Manilius (1807--1869) gave the following
definition of an infinitesimal length in 1850:
\begin{enumerate}\item[]
A length is zero with regard to another, when every finite multiple of
the first is smaller than every finite fraction of the second.%
\footnote{Manilius as translated by Bockstaele\cite[9]{Bo66}.}
\end{enumerate}
Manilius clarifies that ``infinitesimals are not `z\'ero d'une
mani\`ere absolue' but that a finite quantity is never obtained by
taking a finite number of them together, that means that these
quantities are zero with regard to every finite quantity'' (ibid.).
Manilius' definition is essentially a more cumbersome version of the
violation of the Archimedean property.

\subsection{The Leibnizian tradition}
\label{f9}

%Gert Schubring \cite{Sc05} has pursued a different interpretation of
%Cauchy.  Yet in light of the above, it is difficult to accept
%Schubring's claim that at the time ``nobody had thought of the
%calculus in terms of a non-Archimedean continuum.''%
%
%\footnote{\label{f8}Schubring \cite[p.\;2]{Sc22}.  See further in
%  note~\ref{Sc1} on Schubring's position.}
%

No\"el's advocacy of genuine infinitesimals, as described in
Section~\ref{s33b}, was persuasive so long as he was able to remain
faithful to a Leibnizian tradition, such as the distinction between
assignable and inassignable quantities and the definition of
infinitesimal as smaller than any (positive) assignable quantity.

Once No\"el departed from the Leibnizian tradition, his arguments
became less persuasive.  For example, one of the objections by Paque
(see Section~\ref{f29}) concerned a triangle whose base is an
indivisible.  Then the median parallel to the base would have to be
half an indivisible, which Paque saw as a contradiction.  No\"el's
response to this was that there is no contradiction because such a
triangle does not exist.  The pertinent point is not merely that
No\"el's response does not sound convincing to a modern audience, but
rather that No\"el was not faithful to a Leibnizian tradition,
according to which the laws of the finite persist in the infinite,
i.e., the Law of Continuity (LC); see e.g., Katz and Sherry
\cite{13f}.
%Section~\ref{s13}.

On the basis of his LC, Leibniz would have rejected the idea that if
the base is ``too short'' then the triangle does not exist: since for
every ordinary finite base, the triangle does exist, by LC it should
exist for every infinitesimal base, as well.  For such a triangle, the
median still equals half the base, similarly by LC.
%For a discussion of modern formalisations of LC, see
%Section~\ref{s62}.  This is in phi.tex

In fairness to No\"el, it needs to be mentioned that today we have a
clearer understanding of Leibniz's thought based on unpublished
material that was unavailable in mid-19th century Belgium.

\subsection{Lamarle and Paque}
\label{s33}

The opposition to No\"el and Manilius was led by Ernest Lamarle
(1806--1875).  Other critics of infinitesimal methods included
A.\;J.\;N.\;Paque.  Criticisms of infinitesimals generally fall into
two classes: metaphysical and logical.  Berkeley used both (his famous
``ghost'' quip belongs to the metaphysical class).%
\footnote{See Sherry \cite{Sh87}, Katz and Sherry \cite{13f}, and Bair
et al.~\cite[Section~1.1]{22a} for a related analysis of Berkeley's
critique.  See further on Berkeley in note~\ref{f50}.}
Both were used in mid-19th century Belgium.  

Lamarle's and Paque's \emph{metaphysical} criticism was firmly
anchored in pure chimeras (see Section~\ref{s2}).  Thus, Lamarle
wrote:
\begin{enumerate}\item[]
Claiming that the infinitely small constitute special orders of
magnitude amounts, in our view, to reducing them to entities that are
\emph{purely chimerical}.%
\footnote{``Pr\'etendre que les infiniment petits constituent des
ordres particuliers de grandeurs, c'est, selon nous, les r\'eduire \`a
des \^etres purement \emph{chim\'eriques}'' Lamarle \cite[224]{La45};
emphasis added.}
\end{enumerate}
Paque wrote:
\begin{enumerate}\item[]
Formulated in this way and in such terms, the philosophical principle
of the infinitesimal is evidently false; one is therefore within one's
rights to reject in its entirety the set of its mathematical
consequences which do not and cannot possess any existence, quantities
that are truly \emph{purely chimerical}.%
\footnote{``Pos\'e ainsi, tel qu'il l'est et dans de pareils termes,
le principe philosophique infinit\'esimal est \'evidemment faux ; l'on
est donc en droit de repousser de l'ensemble math\'ematique des
auxiliaires qui n'ont et ne peuvent avoir d'existence, de vraies
quantit\'es purement \emph{chim\'eriques}'' Paque
\cite[184--185]{Pa61}; emphasis added.}
\end{enumerate}
For Paque's critique of Leibniz see Section~\ref{f29}.

A \emph{logical} criticism, targeting Leibniz, appeared already on the
second page of Lamarle's lengthy treatise (over 120 pages!), where it
becomes clear that Lamarle is attacking a strawman:
\begin{enumerate}\item[]
The fundamental principle of the method of the infinitely small is the
following:

Two quantities of any order are strictly equal when their difference
is infinitely small of higher order.%
\footnote{``Le principe fondamental de la m\'ethode des infiniment
petits est le suivant:

Deux quantit\'es d'un ordre quelconque sont \emph{rigoureusement
\'egales}, lorsque leur diff\'erence est infiniment petite d'un ordre
sup\'erieur'' Lamarle \cite[222]{La45}; emphasis added.}
\end{enumerate}
Lamarle then proceeds to attack such a notion of strict equality as
``the radical vice represented, from the logical point of view, by the
conception of the infinitely small,''%
\footnote{``le~vice radical que
pr\'esente, au point de vue logique, la conception des infiniment
petits'' (ibid.)}
and poses a rhetorical question:
\begin{enumerate}\item[]
Would it be appropriate that, in a science that is eminently positive
and characterized by an absolute exactness, one should stipulate a
rule whose terms imply a contradiction?%
\footnote{``Convient-il, dans une science \'eminemment positive et
d'une exactitude absolue, de poser une r\`egle dont les termes
impliquent contradiction?''  Lamarle \cite[223]{La45}.}
\end{enumerate}
Accordingly, infinitesimals are not merely chimerical; they are
contradictory.  However, what Lamarle is describing (and attacking) is
not Leibniz's conception.  Leibniz emphasized repeatedly that he
works, not with strict equality implied by Lamarle's wording
\emph{rigoureusement \'egales}, but rather with a generalized relation
of equality ``up to''; see further in Section~\ref{f29}.
%see Section~\ref{s62}, item~\eqref{ione}.  This is in phi.tex

\subsection{Relation of infinite proximity}
\label{f29}

One says that $x$ and~$y$ are in relation of infinite proximity, or
are infinitely close, if~$x-y$ is infinitesimal.  Referring to such a
relation as an ``approximate equality'' has led to numerous
misunderstandings.  Thus, Paque (see Section~\ref{s33}) wrote:
\begin{enumerate}\item[]
Leibniz was forced to respond to accusations which sought no less than
to make out his calculus, lacking a real foundation, to be a veritable
calculus of approximation.%
\footnote{``Leibnitz fut forc\'e de r\'epondre \`a des accusations qui
ne tendaient \`a rien moins qu'\`a faire de son calcul, sans fondement
r\'eel, un v\'eritable calcul d'approximation'' Paque
\cite[180]{Pa61}.}
\end{enumerate}

Paque is clearly using the phrase ``a veritable calculus of
approximation'' in a pejorative sense.  In the same vein, he spoke of
a ``caract\`ere essentiellement approximatif que rev\^et son
hypoth\`ese principale.''%
\footnote{Paque \cite[181]{Pa61}.}
Here Paque conflates the generic meaning of the adjective
``approximate'' and its technical meaning as used to label Leibniz's
relation, more appropriately described as \emph{infinite proximity}.
A relation of infinite proximity can be used rigorously to obtain
precise (and not approximate) results.%
\footnote{Similar misconceptions with regard to Fermat's method of
adequality occur in Breger's work \cite{Br13}; see Katz et
al.~\cite{13e} and Bair et al.~\cite[Section~2.7]{18d}.}

\subsection{Riemann's foundations of infinitesimal geometry}

In his 1854 essay \emph{On the Hypotheses which lie at the Foundations
  of Geometry}, Riemann speaks of the line as being made up of
the~$dx$, describes~$dx$ as `the increments', and speaks of
infinitesimal displacements and of infinitely small quantities such as
\mbox{$x_1 dx_2-x_2\hskip1pt dx_1$}, etc., as well as of infinitely
small triangles.  Analyzing the `next-order term' in the expansion of
the line element on a manifold, Riemann writes:
\begin{enumerate}\item[]
It obviously equals zero if the manifold in question is flat, i.e., if
the square of the line element is reducible to~$\sum dx^2$, and can
therefore be regarded as the measure of deviation from flatness in
this surface direction at this point.  When multiplied by~$-\frac34$
it becomes equal to the quantity which Privy Councilor Gauss has
called the curvature of a surface.%
\footnote{Riemann as translated by Spivak \cite[157]{Sp99}.}
\end{enumerate}
It is apparent from this passage concerning the Gaussian curvature
that Riemann has a pointwise (local) discussion of~$dx$ in mind, and
not merely as an integrand (as in modern notation).

Darrigol notes that some infinitesimal calculations in Riemann's 1861
\emph{Commentatio} memoir were a significant influence on Levi-Civita:
\begin{enumerate}\item[]
Some of Riemann’s formulas involved a peculiar infinitesimal variation
of vectors that Levi-Civita knew to preserve angles and length.  His
interpretation of this variation led him to the concept of parallel
transport.%
\footnote{Darrigol \cite[80]{Da15}.}
\end{enumerate}
Thus the explanatory power of infinitesimals in differential geometry
was already felt by Riemann himself.%
\footnote{How this can be done in modern differential geometry is
explained in Nowik and Katz \cite{15d}.}

\section{Cantor: omega \emph{vs} infinitesimals}

As we saw in Section~\ref{s2b}, Poisson popularized infinitesimals in
France, and both Cauchy and Riemann made routine use of them.  The
prevalent attitudes toward infinitesimals begin to change in the last
quarter of the 19th century.  Bettini's \emph{hallucinations},
Berkeley's \emph{ghosts},%
\footnote{\label{f50}For a recent re-appraisal of Berkeley see
Moriarty's article ``Duelling catechisms: Berkeley trolls Walton on
fluxions and faith'' \cite{Mo23}, as well as \cite{Mo18} and
\cite{Mo22}.}
and d'Alembert--Moigno--Lamarle \emph{chimeras} acquire a novel
identity as \emph{bacilli}%
\footnote{See Ehrlich \cite[4, 55]{Eh06}.}
under the penetrating pen of the great founder of set theory (perhaps
bestowing upon them a superior ontological status to that furnished by
the earlier epithets).

\subsection
{Epithets, from \emph{Abominations} to \emph{Playing with symbols}}
\label{s41}

Cantor used a variety of epithets in reference to infinitesimals.
Already in 1878, Georg Cantor described them as \emph{horribile
dictu}.%
\footnote{See Laugwitz \cite[104]{La02}.}
In a 26\;november\,1890 letter to Veronese, Cantor adheres to the (by
then) entrenched terminology of \emph{ghosts and chimeras}.%
\footnote{See Dauben \cite{Da90}, p.\;235; note\;66, p.\;350.}
But already on 7~september 1890, writing to the same Veronese, he
described them as \emph{abominations}.%
\footnote{Ehrlich \cite[p.\;54 and note 88 there]{Eh06}.}
On 6 october, he accuses Veronese of \emph{playing with symbols}.%
\footnote{\emph{Zeichenspielerei}; cf.~Laugwitz~\cite[120]{La02}.}
By 1893, he pioneers the novel designation of the \emph{infinitary
  Cholera-bacillus}%
\footnote{For the pandemic context of Cantor's remark, see Dauben
\cite{Da90}, p.\;233; note 55, p.\;349.}
%
%and Section~\ref{s51b}. this is in phi.tex
%
in a letter to Vivanti.  In a 1895 letter to Killing, he introduces
the additional epithet \emph{fantasy} (\emph{Phantasterei}).%
\footnote{Dauben \cite{Da90}, note 79, p.\;351.}

Our gentle reader may well wonder wherefore Cantor opposed
infinitesimals.  While he corresponded widely with catholic
theologians,%
\footnote{Tapp \cite{Ta05}.}
there seems to be no indication that Cantor may have been troubled by
eucharist-related issues
% \footnote{See Section~\ref{f1} for related
%theological details.}  This is in phi.tex
that fired the imagination of 17th century opponents such as Jakob
Bidermann.%
\footnote{Festa \cite[198, 207]{Fe92}; Katz et al.~\cite{23g}.}
What, then, was wrong with infinitesimals?

\subsection{Could not have done greater damage}

The answer is laid out clearly in Cantor's 1890 letter to Veronese.
In an earlier letter, Veronese had quoted Otto Stolz to the effect
that Cantor's results have no bearing on Stolz's infinitesimals,
\begin{enumerate}\item[]
because one cannot define their multiplication by the number omega.%
\footnote{Veronese as translated by Ehrlich \cite[54]{Eh06}.}
\end{enumerate}
Stolz and Veronese suggested that there may exist formalisations of
infinity (and infinitesimals) that are different from Cantor's and
admit no multiplication by the Cantorian omega.  Such a suggestion
provoked the following reaction by Cantor on 7 september 1890:
\begin{enumerate}\item[]
Mr.\;Stolz could not have done greater damage to his theory when he
emphasized that his alleged ``\emph{infinitely small magnitudes}''
could not be multiplied by~$\omega$.  An absolute linear magnitude
that is supposed to be real but which nevertheless cannot be
multiplied without limit (and, hence, also with transfinite
multipliers) is not a magnitude at all.%
\footnote{Cantor as translated by Ehrlich \cite[54]{Eh06}.}
\end{enumerate}
The passage suggests that Cantor felt that a theory of infinitesimal
magnitudes that failed to incorporate multiplication by Cantorian
transfinite numbers (e.g.,~$\omega$) in a certain well-behaved
manner,%
\footnote{The relevant conditions ($L_2$ and~$L_3$) were specified in
Ehrlich \cite[46]{Eh06}.}
could not be a valid theory.

One could usefully contrast Cantor's severe verdicts concerning the
work of Stolz and Veronese, with the more liberal account given by
Felix Klein \cite[p.\,156f]{Kl28}.  Klein discusses Veronese's
position as a rejection of the axiom that there is a one-to-one
correspondence between the points of the line and the real numbers.
In a Hilbertian manner, he considers this adoption as legitimate since
it does not lead to contradiction.  But Klein stresses that along with
this kind of admissibility, one also should focus on `fitness to
purpose' (\emph{Zweckm\"a{\ss}igkeit}) when deciding on such
adoptions.  Jahnke and Kr\"omer discuss Klein's usage of
\emph{Zweckm\"a{\ss}igkeit} as a kind of extrinsic justification in
the sense of P.~Maddy; see Jahnke and Kr\"omer \cite{Ja20}.  Klein's
more liberal attitude toward infinitesimals than that of many of his
contemporaries was also analyzed in Bair et al.~\cite{18b}
and~Kanovei et al.~\cite{18i}.%
\footnote{We are grateful to the anonymous referee for providing this
paragraph.}

\subsection{The reverse trajectories of Frege and Peano}  %5

Partly under Cantor's influence, leading logicians, such as Frege and
Peano, wrestled with the question whether infinitesimals are
inconsistent or not, rather than trying to formalize them.  

Frege was initially sympathetic and ended up dismissing them as
inconsistent.  In 1884--85 he did allow for one limited way in which
infinitesimals could be meaningfully defined.  Subsequently in 1903,
he rejected even that approach.%
\footnote{Tappenden \cite{Ta19}.}

Peano, on the contrary, initially endorsed Cantor's deprecation of
infinitesimals, partly due to Peano's rivalry with Veronese%
\footnote{Laugwitz \cite[121--122]{La02}.}
who was Cantor's antagonist.  Eventually Peano himself published an
article in 1910 developing a system with infinitesimals.%
\footnote{Peano \cite{Pe10}; see Bottazzi and Katz \cite{21d} for
details.}
Peano's approach, in modern terms, was to quotient the space of
sequences by the Fr\'echet filter.  Decades later, such an approach
was developed more fully by Schmieden and Laugwitz \cite{Sc58}.

\section{Finale}

Peano's work was one of the first inklings of a re-birth of an
infinitesimalist tradition that had been frowned upon by many
mathematicians since the emergence of the well-known
Cantor--Dedekind--Weierstrass paradigm in analysis during the last
third of the 19th century.  The next major milestones were the work of
Skolem \cite{Sk33}, Hewitt \cite{He48}, {\L}o\'s~\cite{Lo55}, Robinson
\cite{Ro61} and \cite{Ro66}, Hrbacek \cite{Hr78}, Nelson \cite{Ne77}
and others.  Infinitesimal analysis is currently an active field
featuring both mathematical innovations~(see e.g., Hrbacek \cite{24b})
and historical re-appraisal; see e.g., Bair et al.~\cite{20a}, Bair et
al.~\cite{21a}, Katz et al.~\cite{23g}, Bottazzi and Katz \cite{24a},
Ugaglia and Katz \cite{24d}, Katz et al.~\cite{25c}.

\section*{Acknowledgments}

The author is grateful to John T. Baldwin, Karel Hrbacek, Jamie
Tappenden, and the anonymous referee for helpful suggestions.


\begin{thebibliography}{ABCI}


\bibitem[2015]{Al15} Alexander, A, \emph{Infinitesimal--how a
dangerous mathematical theory shaped the modern world}, Scientific
  American/Farrar, Straus and Giroux, New York, 2015.



\bibitem[2024]{Ar24} Arthur, R, and Rabouin, D, `On the unviability of
  interpreting Leibniz's infinitesimals through non-standard
  analysis', \emph{Historia Mathematica} \textbf{66}, (2024), 26--42.





\bibitem[2017a]{17a} Bair, J; B{\l}aszczyk, P; Ely, R; Henry, V;
  Kanovei, V; Katz, K; Katz, M; Kudryk, T; Kutateladze, S; McGaffey,
  T; Mormann, T; Schaps,~D, and Sherry, D, `Cauchy, infinitesimals and
  ghosts of departed quantifiers', \emph{Mat. Stud}.  \textbf{47}
  (2017a), no.\;2, 115--144.
%\url{http://doi.org/10.15330/ms.47.2.115-144},
%\url{https://arxiv.org/abs/1712.00226},
%\url{https://mathscinet.ams.org/mathscinet-getitem?mr=3733080}
  

\bibitem[2017b]{17b} Bair, J; B{\l}aszczyk, P; Ely, R; Henry, V;
  Kanovei, V; Katz, K; Katz, M; Kutateladze, S; McGaffey, T; Reeder,
  P; Schaps, D; Sherry,~D, and Shnider, S, `Interpreting the
  infinitesimal mathematics of Leibniz and Euler', \emph{Journal for
  General Philosophy of Science} \textbf{48} (2017b), no.\;2,
  195--238.
%  \url{http://doi.org/10.1007/s10838-016-9334-z}
%  \url{https://arxiv.org/abs/1605.00455}







\bibitem[2021]{21a} Bair,\;J; B{\l}aszczyk,\;P; Ely,\;R; Katz, M, and
  Kuhlemann,\;K, `Procedures of Leibnizian infinitesimal calculus: An
  account in three modern frameworks', \emph{British Journal for the
  History of Mathematics} \textbf{36} (2021), no.\;3, 170--209.
%  \url{https://doi.org/10.1080/26375451.2020.1851120},
%  \url{https://arxiv.org/abs/2011.12628}
%
%grafting
%


\bibitem[2020a]{20a} Bair, J; B{\l}aszczyk, P; Fuentes Guill\'en, E;
  Heinig, P; Kanovei, V, and Katz, M, `Continuity between Cauchy and
  Bolzano: Issues of antecedents and priority', \emph{British Journal
  for the History of Mathematics} \textbf{35} (2020a), no.\;3,
  207--224.
%https://doi.org/10.1080/26375451.2020.1770015,
%https://arxiv.org/abs/2005.13259

  
\bibitem[2019]{19a} Bair, J; B{\l}aszczyk, P; Heinig, P; Kanovei, V,
  and Katz, M, `19th century real analysis, forward and backward',
  \emph{Antiquitates Mathematicae} \textbf{13} (2019), 19--49.
%  \url{http://doi.org/10.14708/am.v13i1.6440},
%  \url{https://arxiv.org/abs/1907.07451}



\bibitem[2020b]{20b} Bair, J; B{\l}aszczyk, P; Heinig, P; Kanovei, V,
  and Katz, M, `Cauchy's work on integral geometry, centers of
  curvature, and other applications of infinitesimals', \emph{Real
  Analysis Exchange} \textbf{45} (2020b), no.\,1, 127--150.
%  \url{https://doi.org/10.14321/realanalexch.45.1.0127},
%  \url{https://arxiv.org/abs/2003.00438}



\bibitem[2017c]{18b} Bair, J; B{\l}aszczyk, P; Heinig, P; Katz, M;
  Sch\"afermeyer, J, and Sherry,~D, `Klein vs Mehrtens: restoring the
  reputation of a great modern', \emph{Mat. Stud}.  \textbf{48}
  (2017c), no.\;2, 189--219.
%\url{https://arxiv.org/abs/1803.02193},
%\url{http://doi.org/10.15330/ms.48.2.189-219} MR3819950, Zbl
%1459.01017



\bibitem[2022]{22a} Bair, J; Borovik, A; Kanovei, V; Katz, M;
  Kutateladze, S; Sanders, S; Sherry, D, and Ugaglia, M, `Historical
  infinitesimalists and modern historiography of infinitesimals',
  \emph{Antiquitates Mathematicae} \textbf{16} (2022), 189--257.
%  \url{https://arxiv.org/abs/2210.14504},
%  \url{https://doi.org/10.14708/am.v16i1.7169}
%
%reaction
%




\bibitem[2018]{18d} Bair, J; Katz, M, and Sherry, D, `Fermat's
  dilemma: Why did he keep mum on infinitesimals? and the European
  theological context', \emph{Foundations of Science} \textbf{23}
  (2018), no.\;3, 559--595.
%\url{http://doi.org/10.1007/s10699-017-9542-y},
%  \url{https://arxiv.org/abs/1801.00427}


\bibitem[2019]{Ba19} Bair, J, and Mawhin, J, `Le math\'ematicien
  Jean-Nicolas No\"el (1783--1867): un didacticien infinicole du
  19$^e$ si\`ecle', \emph{Revue des Questions scientifiques}
  \textbf{190} (2019), no.\,1--2, 27--59.



\bibitem[2018]{18e} Bascelli, T; B{\l}aszczyk, P; Borovik, A; Kanovei,
  V; Katz, K; Katz, M; Kutateladze, S; McGaffey, T; Schaps, D, and
  Sherry, D, `Cauchy's infinitesimals, his sum theorem, and
  foundational paradigms', \emph{Foundations of Science} \textbf{23}
  (2018), no.\;2, 267--296.
%  \url{http://doi.org/10.1007/s10699-017-9534-y},
%  \url{https://arxiv.org/abs/1704.07723}




\bibitem[1966]{Bo66} Bockstaele, P, `Nineteenth century discussions
  in Belgium on the foundation of the calculus', \emph{Janus}
  \textbf{53} (1966), 1--16.



\bibitem[2021]{21d} Bottazzi, E, and Katz, M, `Infinitesimals via
  Cauchy sequences: refining the classical equivalence', \emph{Open
  Math.}  \textbf{19} (2021), no.\,1, 477--482.



\bibitem[2024]{24a} Bottazzi, E, and Katz, M, `History of Archimedean
  and non-Archimedean approaches to uniform processes: Uniformity,
  symmetry, regularity', \emph{Antiquitates Mathematicae} \textbf{18}
  (2024), 105--148.  \url{https://dx.doi.org/10.14708/am.v18i1.7332},
  \url{https://arxiv.org/abs/2502.18319}



\bibitem[2009]{BS} Bradley, R, and Sandifer, C, \emph{Cauchy's Cours
d'analyse.  An annotated translation}, Sources and Studies in the
  History of Mathematics and Physical Sciences, Springer, New York,
  2009.


\bibitem[2013]{Br13} Breger, H, `Fermat's analysis of extreme values
  and tangents', \emph{Studia Leibnitiana} \textbf{45} (2013), no.\,1,
  20--41.



\bibitem[2019]{Ca19} Cates, D (Tr.), \emph{Cauchy's Calcul
infinit\'esimal, An annotated English translation}, Springer, Cham,
  2019.


\bibitem[1821]{Ca21} Cauchy, A L, \emph{Cours d'Analyse de L'\'Ecole
Royale Polytechnique.  Premi\`ere Partie.  Analyse alg\'ebrique},
  Imprim\'erie Royale, Paris, 1821.\, Also in \emph{Oeuvres
  compl\`etes}, 15 (2.III), Gauthier-Villars, Paris, 1897.



\bibitem[1826]{Ca26} Cauchy, A L, \emph{Le\c cons sur les
applications du calcul infinit\'esimal \`a la g\'eom\'etrie}, Paris,
  Imprim\'erie Royale, 1826.  Also in \emph{Oeuvres Compl\`etes},
  S\'erie 2, Tome\;5.





\bibitem[2000/2004]{Co00} Connes, A, `Cyclic cohomology,
  noncommutative geometry and quantum group symmetries', in Connes et
  al.~\cite{CCG}, pp.\,1--71, 2000 (2004).


\bibitem
[2004] {CCG} Connes, A; Cuntz, J; Guentner, E; Higson, N; Kaminker, J,
and Roberts, J, `Noncommutative geometry, Lectures given at the
C.I.M.E.  Summer School held in Martina Franca, September 3--9, 2000',
edited by S. Doplicher and R.~Longo, \emph{Lecture Notes in
Mathematics}, \textbf{1831}, Springer-Verlag, Berlin; Centro
Internazionale Matematico Estivo (C.I.M.E.), Florence, 2004.



\bibitem[2015]{Da15} Darrigol, O, `The mystery of Riemann's
  curvature', \emph{Historia Mathematica} \textbf{42} (2015), no.\,1,
  47--83.


\bibitem[1990]{Da90} Dauben, J, \emph{Georg Cantor.  His mathematics
and philosophy of the infinite}, Princeton University Press,
  Princeton, NJ, 1990.



\bibitem[2006]{Eh06} Ehrlich, P, `The rise of non-Archimedean
  mathematics and the roots of a misconception. I. The emergence of
  non-Archimedean systems of magnitudes', \emph{Archive for History of
  Exact Sciences} \textbf{60} (2006), no.\,1, 1--121.



\bibitem[1992]{Fe92} Festa, E, `Quelques aspects de la controverse
  sur les indivisibles', in \emph{Geometry and atomism in the Galilean
  school}, 193--207, Bibl. Nuncius Studi Testi, X, Olschki, Florence,
  1992.






\bibitem[2008]{Go08} Goldenbaum, U, and Jesseph, D (eds.),
  \emph{Infinitesimal differences.  Controversies between Leibniz and
  his contemporaries}, De Gruyter, 2008.



\bibitem[1981]{Gr81} Grabiner, J, \emph{The origins of Cauchy's
rigorous calculus}, MIT Press, Cambridge, Mass.--London, 1981.



\bibitem[1990]{Gr90} Grattan-Guinness, I, \emph{Convolutions in
French mathematics, 1800--1840: from the calculus and mechanics to
mathematical analysis and mathematical physics. Vol. II. The turns},
  Science Networks, Historical Studies, 3. Birkh\"auser Verlag, Basel,
  1990.  %MR1119783, Zbl 0836.01013





\bibitem[1948]{He48} Hewitt, E, `Rings of real-valued continuous
  functions.  I', \emph{Transactions of the American Mathematical
  Society} \textbf{64} (1948), 45--99.



\bibitem[1978]{Hr78} Hrbacek, K, `Axiomatic foundations for
  nonstandard analysis', \emph{Fundamenta Mathematicae} \textbf{98}
  (1978), no.\,1, 1--19.



\bibitem[2024]{24b} Hrbacek, K, `Multi-level nonstandard analysis and
  the Axiom of Choice', \emph{Journal of Logic and Analysis}
  \textbf{16}:5 (2024), 1--29.
  %\url{https://doi.org/10.4115/jla.2024.16.5},
%  \url{https://arxiv.org/abs/2405.00621}
Erratum at Hrbacek, K.  `Errata: Multi-level nonstandard analysis and
the axiom of choice', \emph{J. Log. Anal}.  \textbf{16} (2024), Paper
No.\;5c, 3 pp.



\bibitem[2020]{Ja20} Jahnke, H, and Kr\"omer, R, `Rechtfertigen in der
  Mathematik und im Mathematikunterricht', \emph{Journal f\"ur
  Mathematikdidaktik} \textbf{41}(2) (2020), 459--484.



\bibitem[2018]{18i} Kanovei, V; Katz, K; Katz, M, and Mormann, T,
  `What makes a theory of infinitesimals useful? A view by Klein and
  Fraenkel', \emph{Journal of Humanistic Mathematics} \textbf{8}
  (2018), no.\,1, 108--119.
%  \url{http://scholarship.claremont.edu/jhm/vol8/iss1/7},
%  \url{https://arxiv.org/abs/1802.01972}





\bibitem[2013]{13c} Kanovei, V; Katz, M, and Mormann, T, `Tools,
  objects, and chimeras: Connes on the role of hyperreals in
  mathematics', \emph{Foundations of Science} \textbf{18} (2013),
  259--296.
%  \url{http://dx.doi.org/10.1007/s10699-012-9316-5},
%  \url{http://arxiv.org/abs/1211.0244}


\bibitem[2021]{21f} Katz, M, `A two-track tour of Cauchy's
  \emph{Cours}', \emph{Mathematics Today} \textbf{57} (2021), no.\;4,
  154--158.
%  \url{https://arxiv.org/abs/2107.00207}


\bibitem[2025]{23f} Katz, M, and Kuhlemann, K, `Leibniz's contested
  infinitesimals: Further depictions', \emph{Ga\d{n}ita Bh\=arat\=\i}
  (2025), to appear.  \url{https://doi.org/10.32381/GB.2022.45.1.4},
  \url{https://arxiv.org/abs/2501.01193}



\bibitem[2025]{25c} Katz, M; Kuhlemann, K; Sanders, S, and Sherry, D,
  `Formalism~25', \emph{Journal for General Philosophy of Science}
  (2025).  \url{https://arxiv.org/pdf/2502.14811}



\bibitem[2024]{24c} Katz, M; Kuhlemann, K; Sherry, D, and Ugaglia,
  M, `Leibniz on bodies and infinities: \emph{rerum natura} and
  mathematical fictions', \emph{Review of Symbolic Logic} \textbf{17}
  (2024), no.\;1, 36--66.
%\url{https://www.doi.org/10.1017/S1755020321000575},
%\url{https://arxiv.org/abs/2112.08155}
%
%tomasson



\bibitem[2013]{13d} Katz, M, and Leichtnam, E, `Commuting and
  noncommuting infinitesimals', \emph{American Mathematical Monthly}
  \textbf{120} (2013), no.\;7, 631--641.
%  \url{http://doi.org/10.4169/amer.math.monthly.120.07.631},
%  \url{https://arxiv.org/abs/1304.0583}




\bibitem[2013]{13e} Katz, M; Schaps, D, and Shnider, S, `Almost equal:
  The method of adequality from Diophantus to Fermat and beyond',
  \emph{Perspectives on Science} \textbf{21} (2013), no.\;3, 283--324.
%  \url{http://doi.org/10.1162/POSC_a_00101},
%  \url{https://arxiv.org/abs/1210.7750},
%  \url{http://www.ams.org/mathscinet-getitem?mr=3114421}





\bibitem[2013]{13f} Katz, M, and Sherry, D, `Leibniz's infinitesimals:
  Their fictionality, their modern implementations, and their foes
  from Berkeley to Russell and beyond', \emph{Erkenntnis} \textbf{78}
  (2013a), no.\;3, 571--625.
%  \url{http://doi.org/10.1007/s10670-012-9370-y},
%  \url{https://arxiv.org/abs/1205.0174}




\bibitem[2023a]{23g} Katz, M; Sherry, D, and Ugaglia, M, `Of pashas,
  popes, and indivisibles', \emph{Science in Context} \textbf{36}
  (2023a), no.\;2.  \url{https://arxiv.org/abs/2502.11145}



\bibitem[2023b]{23h} Katz, M; Sherry, D, and Ugaglia, M, `When does a
  hyperbola meet its asymptote? Bounded infinities, fictions, and
  contradictions in Leibniz', \emph{Revista Latinoamericana de
  Filosof\'{\i}a} \textbf{49} (2023b), no.\;2, 241--258.
%\url{https://doi.org/10.36446/rlf2023359},
%\url{https://arxiv.org/abs/2311.06023}
%known as delinea




\bibitem[1928]{Kl28} Klein, F, \emph{Elementarmathematik vom h\"oheren
Standpunkt aus III: Pr\"azisions- und Approximationsmathematik},
  Berlin: Springer 3rd ed. 1928.



\bibitem[2019]{La19} Lamand\'e, P, `Sur la conception des objets et
  des m\'ethodes math\'ematiques dans les textes philosophiques de
  d'Alembert', \emph{Historia Mathematica} \textbf{49} (2019), 20--59.


\bibitem[1845]{La45} Lamarle, E, `Essai sur les principes
  fondamenlaux de l'analyse transcendante', \emph{M\'emoires de la
  Soci\'et\'e Royale des Sciences de Li\'ege} \textbf{2} (1845--46),
  221--348.
%acute spelling of Liege is correct according to the pdf sent by Bair




\bibitem[2002]{La02} Laugwitz, D, `Debates about infinity in
  mathematics around 1890: the Cantor--Veronese controversy, its
  origins and its outcome', \emph{NTM (N.S.)}  \textbf{10} (2002),
  no.\;2, 102--126.



\bibitem[1955]{Lo55} {\L}o\'s, J, `Quelques remarques, th\'eor\`emes
  et probl\`emes sur les classes d\'efi\-nissables d'alg\`ebres', in
  \emph{Mathematical interpretation of formal systems}, {98--113},
  North-Holland Publishing, Amsterdam, 1955.
%
%Los, J
%



\bibitem[2017]{Ma17} Mawhin, J, `Fondements du calcul infinitesimal
  en Belgique', Preprint, 12\;dec.\,2017.


\bibitem[2018]{Mo18} Moriarty, C, `The ad hominem argument of
  Berkeley's \emph{Analyst}', \emph{British Journal for the History of
  Philosophy} \textbf{26} (2018), no.\;3, 429--451.


\bibitem[2022]{Mo22} Moriarty, C, `Ructions over fluxions:
  Maclaurin's draft, \emph{The Analyst} Controversy and Berkeley's
  anti-mathematical philosophy', \emph{Studies in History and
  Philosophy of Science} \textbf{96} (2022), 77--86.
%MR4490485


\bibitem[2023]{Mo23} Moriarty, C, `Duelling catechisms: Berkeley
  trolls Walton on fluxions and faith', \emph{Intellectual History
  Review} \textbf{33} (2023), no.\;2, 205--226.
%  \url{https://doi.org/10.1080/17496977.2021.1963933}



\bibitem[1884]{Mo84} \emph{Nature}, L'Abb\'e Moigno (obit.).
  \emph{Nature}, 24 july 1884.



\bibitem[1977]{Ne77} Nelson, E, `Internal set theory: a new approach
  to nonstandard analysis', \emph{Bulletin of the American
  Mathematical Society} \textbf{83} (1977), no.\;6, 1165--1198.


\bibitem[1855]{No55} No\"el, J-N, `Th\'eorie infinit\'esimale
  appliqu\'ee', \emph{M\'emoires de la Soci\'et\'e royale des Sciences
  de Li\'ege. S\'erie 1} \textbf{10} (1855), pp.\;25--136.
%acute spelling of Liege is correct for the period

\bibitem[2015]{15d} Nowik, T, and Katz, M, `Differential geometry
  via infinitesimal displacements', \emph{Journal of Logic and
  Analysis} \textbf{7} (2015), no.\;5, 1--44.
%  \url{http://www.logicandanalysis.com/index.php/jla/article/view/237},
%  \url{https://arxiv.org/abs/1405.0984}



\bibitem[1861]{Pa61} Paque, A J N, `Examen des diverses m\'ethodes
  employ\'ees pour l'\'etablissement et le d\'eveloppement des calculs
  transcendants', \emph{M\'emoires de la Soci\'et\'e Royale des
  Sciences de Li\'ege} \textbf{16} (1861), 145--296.
%
%Alphonse Jean Nicolas Paque according to Bockstaele page 10.
%



\bibitem[1910]{Pe10} Peano, G, `Sugli ordini degli infiniti',
  \emph{Rom. Acc. L. Rend. (5)} \textbf{19} (1910), no.\,1, 778--781.



\bibitem[1961]{Ro61} Robinson, A, `Non-standard analysis',
  \emph{Nederl. Akad. Wetensch. Proc. Ser.\;A} \textbf{64} =
  \emph{Indag. Math.} \textbf{23} (1961), 432--440 (reprinted in
  \emph{Selected Papers}\;\cite{Ro79}, pp.\;3--11).


\bibitem[1966]{Ro66} Robinson, A, \emph{Non-standard analysis},
  North-Holland, Amsterdam, 1966.


\bibitem[1979]{Ro79} Robinson, A, \emph{Selected papers of Abraham
Robinson.  Vol.\;II.  Nonstandard analysis and philosophy}.  Edited
  and with introductions by W. A. J. Luxemburg and S. K\"orner.  Yale
  University Press, New Haven, Conn, 1979.



\bibitem[2018]{18l} Sanders, S, `To be or not to be constructive,
  that is not the question', \emph{Indag. Math. (N.S.)}  \textbf{29}
  (2018), no.\,1, 313--381.
%  \url{https://doi.org/10.1016/j.indag.2017.05.005},
%  \url{https://mathscinet.ams.org/mathscinet-getitem?mr=3739620}



\bibitem[1958]{Sc58} Schmieden, C, and Laugwitz, D, `Eine
  Erweiterung der Infinitesimalrechnung', \emph{Math. Z}.  \textbf{69}
  (1958), 1--39.
  


\bibitem[2005]{Sc05} Schubring, G, \emph{Conflicts between
generalization, rigor, and intuition.  Number concepts underlying the
development of analysis in 17--19th Century France and Germany},
  Sources and Studies in the History of Mathematics and Physical
  Sciences, Springer-Verlag, New York, 2005.



\bibitem[1987]{Sh87} Sherry, D, `The wake of Berkeley's Analyst:
  \emph{rigor mathematicae}?', \emph{Stud. Hist. Philos. Sci.}
  \textbf{18} (1987), no.\;4, 455--480.
%MR0918087, Zbl 0636.01004




\bibitem[1933]{Sk33} Skolem, T, `\"Uber die Unm\"oglichkeit einer
  vollst\"andigen Charakterisierung der Zahlenreihe mittels eines
  endlichen Axiomensystems', \emph{Norsk Mat. Forenings Skr.,
  II. Ser.}, no.\;1/12 (1933), 73--82.


\bibitem[1981]{Sp81} Spalt, D, \emph{Vom Mythos der mathematischen
Vernunft.  Eine Arch\'aologie zum Grundlagenstreit der Analysis oder
Dokumentation einer vergeblichen Suche nach der Einheit der
mathematischen Vernunft}, Wissenschaftliche Buchgesellschaft,
  Darmstadt, 1981.

\bibitem[1999]{Sp99} Spivak, M, \emph{A comprehensive introduction to
differential geometry, Vol.\;II, Third edition}, Publish or Perish,
  Wilmington, Del., 1999.



\bibitem[2005]{Ta05} Tapp, C, \emph{Kardinalit\"at und Kardin\"ale.
Wissenschaftshistorische Aufarbeitung der Korrespondenz zwischen Georg
Cantor und katholischen Theologen seiner Zeit}, Boethius: Texte und
  Abhandlungen zur Geschichte der Mathematik und der
  Naturwissenschaften, 53, Franz Steiner Verlag Wiesbaden GmbH,
  Stuttgart, 2005.


\bibitem[2019]{Ta19} Tappenden, J, `Infinitesimals, magnitudes, and
  definition in Frege', in \emph{Essays on Frege's Basic Laws of
  Arithmetic}, pp.\;235--263, Philip A. Ebert, Marcus Rossberg (Eds.),
  Oxford University Press, 2019.



\bibitem[2024]{24d} Ugaglia, M, and Katz, M, `Evolution of Leibniz's
  Thought in the matter of fictions and infinitesimals', in
  B. Sriraman (ed.), \emph{Handbook of the History and Philosophy of
  Mathematical Practice}, pp.\;341--384, Springer, Cham, 2024.
%\url{https://doi.org/10.1007/978-3-030-19071-2_149-1}
%\url{https://arxiv.org/abs/2310.14249}
%
%phases



\end{thebibliography}
\end{document}